\documentclass[12pt,a4paper,fleqn]{article}
\usepackage{a4wide,amsfonts,amsmath,latexsym,amssymb,euscript,graphicx,units,mathrsfs}

\usepackage{graphicx}
\usepackage{color}
\usepackage{amssymb}
\usepackage{amssymb}
\usepackage[T1]{fontenc}
\usepackage{latexsym}
\usepackage{xypic}
\usepackage{eufrak}
\usepackage{euscript}
\usepackage{amsfonts,amsmath}
\usepackage{verbatim}
\usepackage{fancyhdr}
\usepackage[english]{babel}
\usepackage{mathrsfs}
\usepackage{units}

\newtheorem{prop}{Proposition}[section]
\newtheorem{proposition}{Proposition}[section]

\newtheorem{lemma}[prop]{Lemma}

\newtheorem{theorem}[prop]{Theorem}

\renewcommand{\geq}{\geqslant}
\def\leq{\leqslant}

\newcommand{\R}{\mathbb{R}}

\newcommand{\C}{\mathbb{C}}

\def\1{{\mathbf{1}}}

\def\1{{\mathbf{1}}}
\def\0.5{{\frac{1}{2}}}

\newcommand{\qed}{\nopagebreak\hspace*{\fill}
{\vrule width6pt height6ptdepth0pt}\par}

\newcounter{rea}
\begin{document}

\begin{center}
{\Large{\bf Convergence in law in the second Wiener/Wigner chaos}}
\normalsize
\\~\\ by
Ivan Nourdin\footnote{Email: {\tt inourdin@gmail.com}; IN was partially supported by the
ANR grants ANR-09-BLAN-0114 and ANR-10-BLAN-0121.} and
Guillaume Poly\footnote{Email: {\tt guillaume.poly@crans.org}}\\
{\it Universit\'e de Lorraine and Universit\'e Paris Est}\\~\\
\end{center}

{\small \noindent {{\bf Abstract}:
Let $\mathcal{L}$ be the class of
limiting laws associated with sequences in the second Wiener chaos.
We exhibit a large subset $\mathcal{L}_0\subset\mathcal{L}$ satisfying that,
for any $F_\infty\in\mathcal{L}_0$,
the convergence of only a {\it finite} number of cumulants suffices to
imply the convergence in law of any sequence in the second Wiener chaos to $F_\infty$.
This result is in the spirit of
the seminal paper \cite{nunugio}, in which Nualart and Peccati discovered the surprising fact that
convergence in law for sequences of multiple Wiener-It\^o integrals to the Gaussian is equivalent to convergence of just the fourth cumulant.
Also, we offer analogues of this result
in the case of free Brownian motion and double Wigner integrals, in the context of free probability.
\\

\noindent {\bf Keywords}: Convergence in law; second Wiener chaos; second Wigner chaos; quadratic form; free probability.\\

\noindent {\bf 2000 Mathematics Subject Classification:} 46L54; 60F05, 60G15, 60H05.

\section{Introduction}

Let $X$ be a centered Gaussian process defined on, say, the time interval $[0,1]$.
Consider the sequence of its quadratic variations
\[
V_n=\sum_{k=0}^{n-1} \big(X_{(k+1)/n}-X_{k/n}\big)^2, \quad n\geq 1.
\]
Very often (e.g. when one seeks to build a stochastic calculus with respect to $X$, or
when one wants to estimate an unknown parameter characterizing $X$)
one is interested in determining, provided it exists, the limit in law
of $V_n$ properly normalized, namely,
\begin{equation}\label{fn}
F_n=\frac{V_n-E[V_n]}{\sigma_n},\quad n\geq 1,\quad \mbox{with $\sigma_n= \sqrt{{\rm Var}(V_n)}$}.
\end{equation}
For instance, when $X$ is a fractional Brownian motion of Hurst parameter $H\in(0,1)$, it is known that $F_n$ is asymptotically normal when $H\in(0, 3/4]$, whereas it converges to the Rosenblatt random variable when $H\in(3/4,1)$.
Applying the Gram-Schmidt process leads to the existence of a family $\{e_{k,n}\}_{0\leq k\leq n}$
of $L^2(\R_+)$ satisfying
\[
\int_0^\infty e_{k,n}(x)e_{l,n}(x)dx = E\left[
\big(X_{(k+1)/n}-X_{k/n}\big)
\big(X_{(l+1)/n}-X_{l/n}\big)
\right],\quad 0\leq k,l\leq n-1,
\]
implying in turn that, for any $n\geq 1$,
\[
\big\{X_{(k+1)/n}-X_{k/n}\big\}_{0\leq k\leq n} \overset{\rm law}{=}\left\{
\int_0^\infty e_{k,n}(x)dW_x
\right\}_{0\leq k\leq n},
\]
with $W$ an ordinary Brownian motion.
Hence, since we are concerned with a convergence in law for $F_n$, one can safely replace
its expression (\ref{fn}) by
\[
F_n = \frac1{\sigma_n}\sum_{k=0}^{n-1} \left[
\left(\int_0^\infty e_{k,n}(x)dW_x\right)^2 - \int_0^\infty e_{k,n}(x)^2dx\right].
\]
Next, a straightforward application of the
It\^o's formula yields
\[
F_n = \frac2{\sigma_n}\sum_{k=0}^{n-1}
\int_0^\infty  e_{k,n}(x)dW_x\int_0^x  e_{k,n}(y)dW_y.
\]
That is, $F_n$ has the form of a {\it double integral} with respect to $W$ or, equivalently,
$F_n$ belongs to the {\it second Wiener chaos} associated with $W$.
This simple, illustrating example shows how elements of the second Wiener chaos associated with an
ordinary Brownian motion may be sometimes `hidden' in quantities of interest.

\bigskip

Fix an integer $q\geq 1$. In the late 60's,
Schreiber \cite{Sch} proved that the $q$th Wiener chaos associated with $W$ is
closed under convergence in probability.
In the present paper, we are rather interested in convergence {\it in law}. More
precisely, we aim to answer the following question:
\begin{eqnarray}
\mbox{``{\it Can we describe all the limits in law associated with sequences in the $q$th Wiener chaos?}''}\notag\\
\label{q}
\end{eqnarray}
The answer to (\ref{q}) turns out to be trivial when $q=1$: it is indeed well-known (and easy to prove) that,
 if a sequence of centered Gaussian random variables converge in law, its limit is
centered and Gaussian as well.
In contrast, to provide an answer to (\ref{q}) when $q\geq 3$ is a difficult
problem which is not yet solved.
When $q=2$, the answer to (\ref{q}) is almost contained in reference \cite{russe} by Sevastyanov.
Therein, the author considers quadratic forms of the form
\begin{equation}\label{fd}
F_n= \sum_{i,j=1}^n a_{ij}(N_iN_j-E[N_iN_j]),
\end{equation}
where $N_1,N_2,\ldots\sim\mathcal{N}(0,1)$ are independent,
and characterizes all the possible shapes of any limit in law of a sequence of the form (\ref{fd}).
On the other hand, it
is a classical result that
any element $F$ in the second Wiener chaos can be suitably decomposed as
\begin{equation}\label{deintr}
F\overset{\rm (law)}{=}\sum_{k=1}^\infty \lambda_k (N_k^2-1),
\end{equation}
with $N_1,N_2,\ldots\sim\mathcal{N}(0,1)$ independent and $\{\lambda_k\}_{k\geq 1}\subset\R$
satisfying $\sum_k \lambda_k^2<\infty$.
Thus, sequences in the second Wiener chaos are almost all of the form (\ref{fd})
 and, as a result, the work \cite{russe} is not far away to give a complete answer
to our question (\ref{q}) in the case $q=2$.
In the present paper, we solve (\ref{q}) in full generality for $q=2$, by using another route compared
to \cite{russe}.
Our approach consists to observe that the
identity in law (\ref{deintr}) leads to an expression for the characteristic function $\phi$ of $F$,
 from which it follows that $1/\phi^2$ lies in the
Laguerre-P\'olya class consisting of all real entire functions which are locally uniform limits of real polynomials with real zeros.
Then, relying to a lemma from \cite{HK}, we deduce that the set $\mathcal{L}$ of all limiting laws associated
 to sequences in the second Wiener chaos is composed
of those random variables $F_\infty$ of the form
\begin{equation}\label{formintro}
F_\infty\overset{\rm law}{=}N+G,
\end{equation}
where $G$ lies in the second Wiener chaos
and $N\sim\mathcal{N}(0,\lambda_0^2)$ is independent of $W$.
Thus, for any centered Gaussian process $X$, the limit in law of $F_n$
given by (\ref{fn}) (provided it exists) necessarily has the form (\ref{formintro}).
Another direct consequence of (\ref{formintro}) is that
each random variable in $\mathcal{L}\setminus\{0\}$ has a law which is absolutely
continuous with respect to the Lebesgue measure.

\bigskip

Let us now describe the main result of the present paper.
Relying to (\ref{formintro}), we exhibit
a large subset $\mathcal{L}_0\subset\mathcal{L}$ satisfying that,
for any $F_\infty\in\mathcal{L}_0$,
the convergence of only a {\it finite} number of cumulants suffices to
imply the convergence in law of any sequence in the second Wiener chaos to $F_\infty$.
Theorem \ref{moments-intdouble}, which contains the precise statement, is in the same spirit than
the seminal paper \cite{nunugio}, in which Nualart and Peccati discovered the surprising fact that
convergence in law for sequences of multiple Wiener-It\^o integrals to the Gaussian is equivalent to convergence of just the fourth cumulant.
For an overview of the existing literature around this theme, we refer the reader to the book \cite{NouPecBook}, to the survey \cite{survey}, and to the constantly updated webpage {\tt www.iecn.u-nancy.fr/~nourdin/steinmalliavin.htm}.
\bigskip

Finally, the last part of our paper investigates analogues of all the previous results in the free probability setting.
We are motivated by the fact that there is often a close correspondence between classical
probability and free probability. As one will see, there is no exception to the rule here.
Indeed, let $S=(S_t)_{t\geq 0}$ be a free Brownian motion (defined on a free tracial probability space $(\mathcal{A},\varphi)$) and
let $\mathcal{H}_2^S$ denote the second (symmetric) Wigner chaos associated with $S$, that is,
$\mathcal{H}_2^S$ is the closed linear subspace of $L^2(\varphi)$
generated by the family
\[
\left\{\left(\int_0^\infty h(t)dS_t\right)^2-\int_0^\infty h^2(t)dt: \quad h\in L^2(\R_+)\right\},
\]
where $\int_0^\infty h(t)dS_t$ stands for the Wigner integral of $h$, in the sense of Biane and Speicher \cite{BS}.
Then, if $F_\infty$ is the limit in law
of a given sequence $\{F_n\}$ belonging to $\mathcal{H}^S_2$, we prove that there exists a random variable $A$, distributed according to the semicircular law, as well as
another random variable $G\in \mathcal{H}^S_2$,
 freely independent of $S$, such that
\begin{equation}\label{formintro2}
F_\infty\overset{\rm law}{=}A+G.
\end{equation}
Then, using (\ref{formintro2}) (which is the exact analogue of (\ref{formintro})), we can prove  Theorem \ref{moments-intdouble2}, which is the free counterpart
of Theorem \ref{moments-intdouble}.

\bigskip

The rest of this paper is organized as follows.
In Section 2, we introduce the notation  and we
give several preliminary lemmas which will be used to state and prove our main results in the following sections.
Section 3 is devoted to the proofs of the results described in the introduction in the case
of the standard Brownian motion.
Finally, Section 4 deals with the free probability context.

\section{Preliminaries}

This section gathers the material we shall need to state and prove the
results of this paper.

\subsection{Four objects associated to square-integrable symmetric functions of two variables}

In the whole paper, we shall deal with symmetric functions of $L^2(\R_+^2)$ (that is,
functions $f$ which are square-integrable and that satisfy $f(x,y)=f(y,x)$ a.e.). We denote by
$L^2_s(\R_+^2)$ the set of such functions.

With every function $f\in L^2_s(\R_+^2)$, we associate the selfadjoint Hilbert-Schmidt operator
\[
T_f:L^2(\R_+)\to L^2(\R_+),\quad g\mapsto \int_0^\infty f(\cdot,y)g(y)dy,
\]
as well as the following four objects which are related to it:
\begin{enumerate}
\item[-] We write $\{ \lambda _{k}(f)\}_{k\geq 1}\subset\R$
to indicate the eigenvalues of $T_{f}$ and we assume that $|\lambda_1(f)|\geq |\lambda_2(f)|\geq \ldots$;
\item[-] We set $E(f)=\{\lambda_k(f)\}_{k\geq 1}\setminus\{0\}$;
\item[-] We denote by ${\rm rank}(f)$ the rank of $T_f$, that is, $r={\rm rank}(f)$ if and only if
$\lambda_r(f)\neq 0$ and $\lambda_{r+1}(f)=0$;
\item[-] We denote by $a(f)$ the cardinality of $E(f)$, that is,
$a(f)$ is the number of distinct non-zero eigenvalues of $T_f$.
\end{enumerate}

\subsection{Second Wiener chaos}

Let $W=(W_t)_{t\geq 0}$
be an ordinary Brownian motion on $(\Omega,\mathcal{F},P)$, and
assume that $\mathcal{F}$ is generated by
$W$.
In this subsection, we focus on elements in the second Wiener chaos, that is, we focus on random
variables of the type $F= I_2^W(f)$, with $f\in L^2_s(\R_+^2)$.
For more details on Wiener chaos, we refer to \cite{NouPecBook}.
The following proposition unveils a link between the elements of the second Wiener chaos and some of the objects introduced in the previous subsection.
\begin{proposition}\label{Hi}
For any element $f\in L^2_s(\R_+^2)$, the following equality holds:%
\begin{equation}
I_2^W(f)\overset{\rm (law)}{=}\sum_{k=1}^{\infty }\lambda_{k}(f)\left( N_{k}^{2}-1\right)
\text{,} \label{expX}
\end{equation}%
where $\{N_{k}\}_{k\geq 1}$ is a sequence of
independent $\mathcal{N}(0,1)$ random variables, and the series
converges in ${ L}^{2}(\Omega)$ and almost surely.
\end{proposition}
{\it Proof}. The proof is standard and omitted. See, e.g., \cite[Proposition 2.7.13]{NouPecBook}.\qed

\bigskip

By using (\ref{expX}), it is easily seen that
the cumulants of $F$ are given by
\[
\kappa_1(I_2^W(f))=0\quad\mbox{and}\quad \kappa_r(I_2^W(f))=2^{r-1}(r-1)! \sum_{k=1}^\infty \lambda_k(f)^r,\quad r\geq 2.
\]
See, e.g., \cite[Identities (2.7.17)]{NouPecBook}.

\bigskip

We shall also need the following hypercontractivity property.
\begin{theorem}\label{hyper}
Let $F$ be a double Wiener-It\^o integral. Then, for all $r\geq 1$, we have
\begin{equation}\label{hyper-e}
E[|F|^r]\leq (r-1)^r E[F^2]^{r/2}.
\end{equation}
\end{theorem}
{\it Proof}. See, e.g., \cite[Corollary 2.8.14]{NouPecBook}.
\qed

\subsection{Second Wigner chaos}

Let $S=(S_t)_{t\geq 0}$ be a free Brownian motion, defined on a non-commutative probability space
$(\mathcal{A},\varphi)$. That is, $S$ is a stochastic process
starting from 0,
with freely independent increments, and such that $S(t)-S(s)\sim\mathcal{S}(0,t-s)$ is
a centered semicircular random variable with variance $t-s$ for all $t\geq s$.
(We may think of free Brownian motion as `infinite-dimensional matrix-valued Brownian motion'.)
For more details about the construction and features of $S$, see \cite[Section 1.1]{BS}
and the references therein.

In this subsection, we focus on elements in the (symmetric) second Wigner chaos, that is, we focus on random
variables of the type $F= I_2^S(f)$, with $f\in L^2_s(\R_+^2)$.
The following result is nothing but the free counterpart of Proposition \ref{Hi}.

\begin{proposition}\label{Hi2}
For any element $f\in L^2_s(\R_+^2)$, the following equality in law holds:
\begin{equation}
I_2^S(f)\overset{\rm law}{=}\sum_{k=1}^{\infty }\lambda_{k}(f)\left( S_{k}^{2}-1\right)
\text{,} \label{expX2}
\end{equation}%
where $\{S_{k}\}_{k\geq 1}$ is a sequence of
freely independent
centered semicircular random variables with unit variance, and the series
converges in ${ L}^{2}(\varphi)$.
\end{proposition}
{\it Proof}. It is an immediate extension of \cite[Proposition 2.7.13]{NouPecBook}
to the free case.\qed

\bigskip

An immediate consequence of Proposition \ref{Hi2} is that the free cumulants (see \cite{NS}) of $F$  are given by
\[
\widehat{\kappa}_1(I_2^S(f))=0\quad\mbox{and}\quad \widehat{\kappa}_r(I_2^S(f))=
\sum_{k=1}^\infty \lambda_k(f)^r,\quad r\geq 2.
\]

\section{Our results in the classical Brownian motion case}

Let $W=(W_t)_{t\geq 0}$
be an ordinary Brownian motion on $(\Omega,\mathcal{F},P)$.
The next theorem describes all the limits in law associated with sequences in the second Wiener chaos.

\begin{theorem}\label{structure}
Let $\{F_n\}_{n\geq 1}$ be a sequence of double Wiener integrals that converges in law to $F_\infty$.
Then, there exists $\lambda_0\in\R$ and $f\in L^2_s(\R_+^2)$ such that
\[
F_\infty \overset{\rm (law)}{=}N + I_2^W(f),
\]
where $N\sim\mathcal{N}(0,\lambda_0^2)$ is independent of the underlying Brownian motion $W$.
\end{theorem}

 During the proof of Theorem \ref{structure}, we shall need the following result taken from \cite{HK} (more precisely, it is a suitable combination of Lemma 1 and Lemma 2 in \cite{HK}).

\begin{lemma}\label{LPlemma}
Let $\{G_n\}$ be a sequence of entire functions of
the form
\[
G_n(z)=e^{\alpha_n z+\beta_n}\prod_{k=1}^\infty \left(1-z/z_{k,n}\right)e^{z/z_{k,n}},\quad z\in\C,
\]
with $\alpha_n,\beta_n\in\C$, and where the zeros $\{z_{k,n}\}_{k\geq 1}$ of $G_n$ are included in $\R\setminus\{0\}$ and satisfy the condition
\[
\sum_{k=1}^\infty |z_{k,n}|^{-2}\leq M,\quad n=1,2,\ldots,
\]
for some constant $M$ independent of $n$.
Assume that $G_n$ converges uniformly on a disc about the origin, to a limit function
$\not\equiv 0$. Then $G_n$ converge uniformly on every bounded set, to an entire function $G_\infty$
of the form
\begin{equation}\label{form}
G_\infty(z)=e^{az^2+bz+c}\prod_{k=1}^\infty \left(1-z/z_k\right)e^{z/z_k},
\end{equation}
where $a,b,c\in\C$ and where the $z_k$ are real and such that $\sum_{k=1}^\infty |z_k|^{-2}<\infty$.
\end{lemma}

We shall also need the following result, which is a straightforward consequence of the Paley inequality
(\ref{papalele}) as well as the hypercontractivity property (\ref{hyper-e}).
\begin{lemma}\label{paley}
Let $\{F_n\}_{n\geq 1}$ be a tight sequence of double Wiener-It\^o integrals.
Then
\[
\sup_{n\geq 1}E[|F_n|^p]<\infty
\quad
\mbox{for all $p\geq 1$}.
\]
\end{lemma}
{\it Proof}. Let $Z$ be a positive random variable such that $E[Z]=1$ and let $\theta\in(0,1)$.
Consider the decomposition $Z=Z{\bf 1}_{\{Z>\theta\}}+Z{\bf 1}_{\{Z\leq\theta\}}$ and
take the expectation. One deduces, using Cauchy-Schwarz, that
$\sqrt{E[Z^2]}\sqrt{P(Z>\theta)}+\theta\geq 1$,
that is (Paley inequality),
\begin{equation}\label{papalele}
E[Z^2]\,P(Z>\theta)\geq (1-\theta)^2.
\end{equation}
On the other hand, we have by hypercontractivity (\ref{hyper-e}) that
$E\big[F_n^4\big]\leq 81\,E\big[F_n^2\big]^2$.
Combining this latter fact with (\ref{papalele})
yields that, for all $\theta \in (0,1)$ and with $Z=F_n^2/E[F_n^2]$,
\begin{equation}\label{Paley}
P\big( F_n^2 > \theta E [F_n^2] \big) \geq \frac1{81}(1-\theta)^2 .
\end{equation}
The sequence $\{F_n\}_{n\geq 1}$ being tight, one can choose $M>0$
 large enough so that $P(F_n^2 > M) \leq \frac{1}{324}$
for all $n\geq 1$.
By applying \eqref{Paley} with $\theta=1/2$, one gets that
\begin{align*}
P\big(F_n^2 > M\big) \leq \frac{1}{324} \leq P\big( F_n^2 > \frac12\, E [F_n^2] \big),
\end{align*}
from which one deduce that $\sup_{n\geq 1} E [F_n^2] \leq 2M<\infty$.
The desired conclusion follows (once again!) from the hypercontractivity property (\ref{hyper-e}).
\qed

\bigskip

We are now in a position to prove Theorem \ref{structure}.

\bigskip

\noindent
{\it Proof of Theorem \ref{structure}}.
We first observe that $\sup_{n\geq 1}E[|F_n|^p]<\infty$ for all $p\geq 1$ by Lemma \ref{paley}. Hence, without loss of generality, we may and do assume that $\sup_{n\geq 1}E[F_n^2]\leq 1$.

The rest of the proof is divided into two steps.

\bigskip
{\it Step 1}. We claim that $\phi_n$, defined on $\Omega=\{z\in\C:\,|{\rm Re}z|< e^{-1}\}$ as
\[
\phi_n(z)=E[e^{zF_n}],
\]
 converges uniformly on compact sets of $\Omega$, towards an holomorphic
function noted $\phi_\infty$.
First, observe that the function $\phi_n$ is well-defined and uniformly bounded on compact sets of $\Omega$. Indeed, by hypercontractivity (\ref{hyper-e})
and for every $q>2$,
\[
E[|F_n|^q]^{1/q} \leq q-1.
\]
Thus, $P[|F_n|> u]\leq u^{-q} (q-1)^q$ for every $u>0$. Choosing $q = q(u) = 1+u/e$, the previous relation shows that $P[|F_n|> u]\leq e^{-u/e}$ for every $u> e$. By a Fubini argument,
\begin{eqnarray}
\sup_{n\geq 1}E[|e^{zF_n}|]\leq \sup_{n\geq 1}E[e^{|{\rm Re}z|\,|F_n|}] =\sup_{n\geq 1}\left(1+ |{\rm Re}z| \int_0^\infty e^{|{\rm Re}z|u} P\big[| F_n | >u \big] du\right)<\infty\label{laplace}
\end{eqnarray}
for any $z\in\C$ such that $|{\rm Re}z|<e^{-1}$.
Moreover, since
\[
\int_{\partial T}\phi_n(z)dz = \int_{\partial T}E[e^{zF_n}]dz = E\left[\int_{\partial T}e^{zF_n}dz\right]=0
\]
for all triangle $\partial T\subset\Omega$, one deduces from the Goursat theorem that
each $\phi_n$ is holomorphic on $\Omega$.

Since the sequence $\{\phi_n\}_{n\geq 1}$
is bounded on compact sets of $\Omega$, one can assume by Montel's theorem that it converges uniformly on
compact sets of $\Omega$ (towards, say, $h$) and one is thus left to show that the limit is unique.
Since $F_n\overset{\rm law}{\to}F_\infty$, we have that $\phi_n(it)$ converges pointwise to $\phi_\infty(it)$
for all $t\in\R$
(with obvious notation). We deduce that $h(it)=\phi_\infty(it)$ for all $t\in\R$.
But, as an immediate consequence of (\ref{laplace}) and of the continuous mapping theorem, we have that the function $\phi_\infty(z)=E[e^{zF_\infty}]$ is well-defined and holomorphic on $\Omega$. Hence, $h=\phi$ on $\Omega$ and the desired claim follows.

\bigskip

{\it Step 2}. Let $\{N_k\}_{k\geq 1}$ denote a sequence of independent $\mathcal{N}(0,1)$ random variables.
Since $F_n=I_2(f_n)$ is a double Wiener-It\^o integral, its law can be expressed (see Proposition \ref{Hi}) as
\[
F_n\overset{\rm law}{=} \sum_{k=1}^\infty \lambda_{k,n}(N_k^2-1)
\]
where the $\lambda_{k,n}$ are the eigenvalues of the Hilbert-Schmidt operator
$T_n:L^2(\R_+)\to L^2(\R_+)$
    defined as $T_n(g)=\langle f_n,g\rangle_{L^2(\R_+)}$.
    In particular, $\sum_{k=1}^\infty \lambda_{k,n}^2 = \frac12E[F_n^2]<\infty$.
It is straightforward to check that
\[
G_n(z):=\frac{1}{\phi_n(z)^2}=\prod_{k=1}^\infty \left(1-2\lambda_{k,n}z\right)e^{2\lambda_{k,n}z},\quad z\in\Omega.
\]
Since $\phi_\infty$ is not identically zero and $\phi_n(z)\neq 0$ for all $z\in\Omega$,
the Hurwitz principle applies and yields that $\phi_\infty(z)\neq 0$ for all $z\in\Omega$.
Therefore,
$G_n$ converges uniformly on compact sets of $\Omega$ to $G_\infty=\phi_\infty^{-2}$.
Thanks to Lemma \ref{LPlemma}, we deduce that $G_\infty$ has the form
\begin{equation}\label{form21}
G_\infty(z)=e^{az^2+bz+c}\prod_{k=1}^\infty \left(1-2\lambda_{k}z\right)e^{2\lambda_{k}z},\quad z\in\Omega,
\end{equation}
where $a,b,c\in\C$ and where the $\lambda_k$ are real and such that $\sum_{k=1}^\infty \lambda_k^2<\infty$.
In (\ref{form21}), we necessarily have $a\in\R_+$, $b=0$ and $e^c=1$.
Indeed, $e^c=G_\infty(0)=1$. Moreover,
\[
\prod_{k=1}^\infty \left(1-2\lambda_{k}z\right)e^{2\lambda_{k}z} = e^{
-\sum_{j=2}^\infty \frac{(2z)^j}{j}\sum_{k=1}^\infty \lambda_k^j},\quad z\in\Omega,
\]
so that
\[
b=G'_\infty(0)=-2E[F_\infty]=0
\]
and
\[
2a-4\sum_{k=1}^\infty\lambda_k^2 = G''_\infty(0)=-2E[F_\infty^2],
\]
implying that $a=2\sum_{k=1}^\infty \lambda_k^2 - E[F_\infty^2]\in\R$.
Moreover, using Fatou's lemma we deduce that
\[
2\sum_{k=1}^\infty \lambda_k^2 \leq 2\liminf_{n\to\infty} \sum_{k=1}^\infty \lambda_{m^{k}_{n},n}^2
\leq 2\liminf_{n\to\infty} \sum_{k=1}^\infty \lambda_{k,n}^2
= \liminf_{n\to\infty} E[F_n^2]=E[F_\infty^2].
\]
That is, $a\leq 0$.
Consequently, (\ref{form21}) implies that
\begin{equation}\label{form2}
F_\infty \overset{\rm law}{=} N + \sum_{k=1}^\infty \lambda_k (N_k^2-1),
\end{equation}
where $N\sim\mathcal{N}(0,-a)$
is independent from $N_k$, $k\geq 1$.
To conclude the proof, it suffices to observe that
\[
\sum_{k=1}^\infty \lambda_k (N_k^2-1) \overset{\rm law}{=} I_2^W(f),
\]
where $f$ is given by $f(x,y)=\sum_{k=1}^\infty e_k(x)e_k(y)$, with $\{e_k\}_{k\geq 1}$
any orthonormal basis of $L^2(\R_+)$.

\qed

\bigskip

In the next theorem, we exhibit
a large subset $\mathcal{L}_0\subset\mathcal{L}$ such that,
for any $F_\infty\in\mathcal{L}_0$,
the convergence of only a {\it finite} number of cumulants suffices to
imply the convergence in law of any sequence in the second Wiener chaos to $F_\infty$.

\begin{theorem}\label{moments-intdouble}
Let $f\in L^2_s(\R_+^2)$ with $0\leq {\rm rank}(f)<\infty$,
let $\mu_0\in\R$ and let $N\sim\mathcal{N}(0,\mu_0^2)$ be independent of the underlying Brownian motion $W$.
Assume that $|\mu_0|+\|f\|_{L^2(\R_+)}>0$ and set
\[
Q(x)=x^{2(1+{\bf 1}_{\{\mu_0\neq 0\}})}\prod_{i=1}^{a(f)}(x-\lambda_i(f))^2.
\]
Let $\{F_n\}_{n\geq 1}$ be a sequence of double Wiener-It\^o integrals.
Then, as $n\to\infty$, we have 
\begin{enumerate}
\item[(i)] $F_n\overset{\rm law}{\to} N+I_2^W(f)$
\end{enumerate}
if and only if all the following are satisfied:
\begin{enumerate}
\item[(ii-a)] $\kappa_2(F_n)\to \kappa_2(N+I_2^W(f))=\mu_0^2+2\|f\|^2_{L^2(\R_+^2)}$;
\item[(ii-b)] $\sum_{r=3}^{{\rm deg}Q} \frac{Q^{(r)}(0)}{r!}\,\frac{\kappa_r(F_n)}{(r-1)!2^{r-1}}
\to \sum_{r=3}^{{\rm deg}Q} \frac{Q^{(r)}(0)}{r!}\,\frac{\kappa_r(I_2^W(f))}{(r-1)!2^{r-1}}$;
\item[(ii-c)] $\kappa_r(F_n)\to \kappa_r(I_2^W(f))$ for  $a(f)$ consecutive values of $r$, with $r\geq 2(1+{\bf 1}_{\{\mu_0\neq 0\}})$.
\end{enumerate}
\end{theorem}

Before doing the proof of Theorem \ref{moments-intdouble}, let us detail two explicit examples.
\begin{enumerate}
\item
Consider first the situation where $\|f\|_{L^2(\R_+^2)}=0$ (that is, ${\rm rank}(f)=a(f)=0$) and $\mu_0\neq 0$. In this case, $Q(x)=x^4$ and condition $(ii-c)$ is
immaterial. Therefore, $F_n\overset{\rm law}{\to}
\mathcal{N}(0,\mu_0^2)$ if and only if $\kappa_2(F_n)\to \mu_0^2$ (condition $(ii-a)$) and $\kappa_4(F_n)\to 0$ (conditions $(ii-b)$). As such, one recovers the celebrated Nualart-Peccati criterion (in the case of double integrals), see \cite{nunugio}.
\item
Consider now the situation where, in Theorem \ref{moments-intdouble}, one has $r={\rm rank}(f)<\infty$, $a(f)=1$, $\lambda_1(f)=1$ and $\mu_0=0$. This corresponds to the case where the limit $I^W_2(f)$
is a centered chi-square random variable with $r$ degrees of freedom. In this case, $Q(x)=x^4-2x^3+x^2$.
Therefore,  $F_n\overset{\rm law}{\to}
I^W_2(f)$ if and only if $\kappa_2(F_n)\to 2r$ (conditions $(ii-a)$ and $(ii-c)$) and
$\kappa_4(F_n)-12\,\kappa_3(F_n)\to -48r$ (condition $(ii-b)$). As such, one recovers a theorem of Nourdin and Peccati in the special case of double integrals, see \cite{noncentral}.
\end{enumerate}

For the proof of Theorem \ref{moments-intdouble}, we shall need the following auxiliary lemma.

\begin{lemma}\label{vandermonde}
Let $\mu_0\in\R$, let $a\in\mathbb{N}^*$, let $\mu_1,\ldots,\mu_a\neq 0$ be pairwise distinct real numbers,  and let
$m_1,\ldots,m_a\in\mathbb{N}^*$.
Set \[
Q(x)=x^{2(1+{\bf 1}_{\{\mu_0\neq 0\}})}\prod_{i=1}^a(x-\mu_i)^2
\]
Assume that $\{\lambda_j\}_{j\geq 0}$ is a square-integrable sequence of real numbers satisfying
\begin{eqnarray}
\label{eq1}
\lambda_0^2+\sum_{j=1}^\infty \lambda_j^2&=&\mu_0^2+\sum_{i=1}^a m_i\,\mu_i^2\\
\label{eq3}
\sum_{r=3}^{2(1+{\bf 1}_{\{\mu_0\neq 0\}}+a)}\frac{Q^{(r)}(0)}{r!}\sum_{j=1}^\infty \lambda_j^r
&=&
\sum_{r=3}^{2(1+{\bf 1}_{\{\mu_0\neq 0\}}+a)}\frac{Q^{(r)}(0)}{r!}\sum_{i=1}^a m_i\,\mu_i^r\\
\sum_{j=1}^\infty \lambda_j^r&=&\sum_{i=1}^a m_i\,\mu_i^r\quad \mbox{for `$a$' consecutive values of $r\geq 2(1+{\bf 1}_{\{\mu_0\neq 0\}})$}.\notag\\
\label{eq2}
\end{eqnarray}
Then:
\begin{itemize}
\item[(i)] $|\lambda_0|=|\mu_0|$.
\item[(ii)] The cardinality of the set $S=\{j\geq 1:\,\lambda_j\neq 0\}$ is finite.
\item[(iii)] $\{\lambda_j\}_{j\in S}=\{\mu_i\}_{1\leq i\leq a}$.
\item[(iv)] for any $i=1,\ldots,a$, one has $m_i=\#\{j\in S:\,\lambda_j=\mu_i\}$.
\end{itemize}
\end{lemma}
{\it Proof}.
We divide the proof according to the nullity of $\mu_0$.

\bigskip

{\it First case}: $\mu_0=0$.
We have $Q(x)=x^{2}\prod_{i=1}^a(x-\mu_i)^2$.
Since the polynomial $Q$ can be rewritten as
\[
Q(x)=\sum_{r=2}^{2(1+a)}\frac{Q^{(r)}(0)}{r!}x^r,
\]
assumptions (\ref{eq1}) and (\ref{eq2}) together ensure that
\[
\lambda_0^2\prod_{i=1}^a\mu_i^2+\sum_{j=1}^\infty Q(\lambda_j)=\sum_{i=1}^a m_i Q(\mu_i)=0.
\]
Because $Q$ is positive and $\prod_{i=1}^a\mu_i^2\neq 0$, we deduce that $\lambda_0=0$ and $Q(\lambda_j)=0$ for all $j\geq 1$, that is,
$\lambda_j\in\{0,\mu_1,\ldots,\mu_a\}$ for all $j\geq 1$. This shows claims $(i)$ as well as:
\begin{equation}\label{inclusion}
\{\lambda_j\}_{j\in S}\subset\{\mu_i\}_{1\leq i\leq a}.
\end{equation}
Moreover, since the sequence $\{\lambda_j\}_{j\geq 1}$ is square-integrable, claim $(ii)$ holds true as well.
It remains to show $(iii)$ and $(iv)$. For any $i=1,\ldots,a$, let $n_i=\#\{j\in S:\,\lambda_j=\mu_i\}$.
Also, let $r\geq2$ be such that $r,r+1,\ldots,r+a-1$ are `$a$' consecutive values satisfying (\ref{eq3}).
We then have
\[
\left(
\begin{array}{cccc}
\mu_1^{r}&\mu_2^{r}&\cdots&\mu_a^{r}\\
\mu_1^{r+1}&\mu_2^{r+1}&\cdots&\mu_a^{r+1}\\
\vdots&\vdots&\ddots&\vdots\\
\mu_1^{r+a-1}&\mu_2^{r+a-1}&\cdots&\mu_a^{r+a-1}
\end{array}
\right)
\left(\begin{array}{c}
n_1-m_1\\
n_2-m_2\\
\vdots\\
n_a-m_a
\end{array}
\right)
=
\left(\begin{array}{c}
0\\
0\\
\vdots\\
0
\end{array}
\right).
\]
Since $\mu_1,\ldots,\mu_a\neq 0$ are pairwise distinct, one has (Vandermonde matrix)
\begin{eqnarray*}
&&\det\left(
\begin{array}{cccc}
\mu_1^{r}&\mu_2^{r}&\cdots&\mu_a^{r}\\
\mu_1^{r+1}&\mu_2^{r+1}&\cdots&\mu_a^{r+1}\\
\vdots&\vdots&\ddots&\vdots\\
\mu_1^{r+a-1}&\mu_2^{r+a-1}&\cdots&\mu_a^{r+a-1}
\end{array}
\right)\\
&=&
\prod_{i=1}^a \mu_i^r
\times
\det\left(
\begin{array}{cccc}
1&1&\cdots&1\\
\mu_1&\mu_2&\cdots&\mu_a\\
\vdots&\vdots&\ddots&\vdots\\
\mu_1^{a-1}&\mu_2^{a-1}&\cdots&\mu_a^{a-1}
\end{array}
\right)
\neq 0,
\end{eqnarray*}
from which $(iv)$ follows. Finally, recalling the inclusion (\ref{inclusion}) we deduce $(iii)$.

\bigskip

{\it Second case}: $\mu_0\neq 0$.
In this case, one has $Q(x)=x^{4}\prod_{i=1}^a(x-\mu_i)^2$
and claims $(ii)$, $(iii)$ and $(iv)$ may be shown by following the same line of reasoning as above. We then deduce claim $(i)$ by looking at (\ref{eq1}).
\qed

\bigskip

We are now in a position to prove Theorem \ref{moments-intdouble}.

\bigskip

\noindent
{\it Proof of Theorem \ref{moments-intdouble}}.
The implications $(i)\to(ii-a),(ii-b),(ii-c)$ are immediate consequences of the Continuous Mapping Theorem together with Lemma \ref{paley}.
Now, assume $(ii-a),(ii-b),(ii-c)$ and let us show that $(i)$ holds true.
The sequence $\{F_n\}$ being bounded in $L^2$, it is tight by Prokhorov's theorem.
 Hence, to prove claim $(i)$ it is sufficient to show that any subsequence $\{F_{n'}\}$ converging in law to some
 random variable $F_\infty$ is necessarily such that $F_\infty\overset{\rm law}{=}N+I_2^W(f)$.
 From now on, and only for notational convenience, we assume that $\{F_n\}$ itself converges to $F_\infty$.
 By hypercontractivity (\ref{hyper-e}), one has that $\kappa_r(F_n)\to \kappa_r(F_\infty)$ for all $r$.
  Thanks to Theorem \ref{structure}, we know that
  \begin{equation}\label{nec}
  F_\infty \overset{\rm law}{=} \lambda_0 N_0 + \sum_{j=1}^\infty \lambda_j (N_j^2-1),
  \end{equation}
  for some $\lambda_0,\lambda_1,\ldots$ satisfying $\sum_{j}\lambda_j^2<\infty$ and where $N_0,N_1,N_2,\ldots\sim\mathcal{N}(0,1)$
  are independent. Combining our assumptions $(ii-a)$,$(ii-b)$,$(ii-c)$ with (\ref{nec}), we deduce that (\ref{eq1})-(\ref{eq3})-(\ref{eq2}) hold true.
  If $a\neq 0$, we deduce from Lemma \ref{vandermonde} that
 $G_\infty \overset{\rm law}{=}F_\infty$, thereby concluding the proof of Theorem \ref{moments-intdouble}.
 If $a=0$ and $\mu\neq 0$, then the desired conclusion follows from the Nualart-Peccati criterion of asymptotic normality.
 Finally, if $a=0$ and $\mu=0$, then (\ref{eq1}) implies that $F_\infty=0$ a.s., and the proof of Theorem \ref{moments-intdouble} is also concluded in this case.
\qed

\section{Our results in the free Brownian motion case}

Let $S=(S_t)_{t\geq 0}$ be a free Brownian motion, defined on a non-commutative probability space
$(\mathcal{A},\varphi)$.
The following result fully describes all the possible limiting laws for sequences in the second (symmetric) Wigner chaos.
It is the exact analogue of Theorem \ref{structure}.

\begin{theorem}\label{structure2}
Let $\{F_n\}_{n\geq 1}$ be a sequence in the second Wigner chaos that converges in law to $F_\infty$.
Then, there exists $\lambda_0\in\R$ and $f\in L^2_s(\R_+^2)$ such that
\[
F_\infty \overset{\rm law}{=}A + I_2^S(f),
\]
where $A\sim\mathcal{S}(0,\lambda_0^2)$ is freely independent of the underlying free Brownian motion $S$.
\end{theorem}
{\it Proof}.
We have $F_n=I_2^S(f_n)$ with $f_n\in L^2_s(\R_+^2)$ for each $n\geq 1$.
Set $G_n=I_2^W(f_n)$ where $W=(W_t)_{t\geq 0}$ stands for an ordinary Brownian motion. We have 
\[
E[G_n^2]=2\|f_n\|^2_{L^2(\R_+^2)}=2\varphi(F_n^2)
\to 2\varphi(F_\infty^2)\] 
as $n\to\infty$,
so $\{G_n\}$ is bounded in $L^2(\Omega)$.
By Prokhorov's theorem, it is tight and there exists a subsequence $\{G_{n'}\}$
that converges in law to, say, $G_\infty$. By Theorem \ref{structure}, the limit
$G_\infty$ has necessarily the form $G_\infty= N + I_2^W(f)$ with $f\in L^2(\R_+^2)$
and $N\sim\mathcal{N}(0,2\lambda_0^2)$ independent of $W$. Using Lemma \ref{paley},
we deduce that, as $n'\to\infty$,
\[
2\sum_{k=1}^\infty \lambda_k(f_{n'})^2=\kappa_2(G_{n'})\to \kappa_2(G_\infty)=
2\lambda_0^2+2\sum_{k=1}^\infty \lambda_k(f)^2
\]
and
\[
2^{r-1}(r-1)!\sum_{k=1}^\infty \lambda_k(f_{n'})^r=\kappa_r(G_{n'})\to \kappa_r(G_\infty)=
2^{r-1}(r-1)!\sum_{k=1}^\infty \lambda_k(f)^r,\quad r\geq 3.
\]
Now, since $F_n$ converges in law to $F_\infty$, we have, as $n\to\infty$,
\[
\sum_{k=1}^\infty \lambda_k(f_{n})^r=\widehat{\kappa}_r(F_n)\to \widehat{\kappa}_r(F_\infty)\quad\mbox{for any $r\geq 2$}.
\]
We deduce that
\[
\widehat{\kappa}_2(F_\infty)=\lambda_0^2+\sum_{k=1}^\infty\lambda_k(f)^2\quad\mbox{and}\quad
\widehat{\kappa}_r(F_\infty)=\sum_{k=1}^\infty\lambda_k(f)^r,\quad
r\geq 3,
\]
from which the announced claim follows.
\qed

\bigskip

Using cumulants (as in the previous proof, see also Proposition \ref{paley}), it is straightforward to check the following
Wiener-Wigner transfer principle between limits in the second Wiener chaos and limits in the second Wigner chaos.

\begin{theorem}\label{transfer}
Let $\{f_n\}_{n\geq 1}$ be a sequence of elements of $L^2_s(\R_+^2)$,
let $W=(W_t)_{t\geq 0}$ be an ordinary Brownian motion and let $S=(S_t)_{t\geq 0}$ be a free
Brownian motion.
Let $N\sim\mathcal{N}(0,2\lambda_0^2)$
(resp. $A\sim\mathcal{S}(0,\lambda_0^2)$) be independent of $W$ (resp. $S$).
Then, as $n\to\infty$,  $I_2^W(f_n)\overset{\rm law}{\to} N+I_2^W(f)$
if and only if  $I_2^S(f_n)\overset{\rm law}{\to} A+I_2^S(f)$.
\end{theorem}

Also, reasoning as in the proof of Theorem \ref{moments-intdouble},
we may obtain the following characterization in terms of a finite number of cumulants
for convergence in law within the second Wigner chaos.

\begin{theorem}\label{moments-intdouble2}
Let $f\in L^2_s(\R_+^2)$ with $0\leq {\rm rank}(f)<\infty$,
let $\mu_0\in\R$ and let $A\sim\mathcal{S}(0,\mu_0^2)$ be independent of the underlying free Brownian motion $S$.
Assume that $|\mu_0|+\|f\|_{L^2(\R_+)}>0$ and set
\[
Q(x)=x^{2(1+{\bf 1}_{\{\mu_0\neq 0\}})}\prod_{i=1}^{a(f)}(x-\lambda_i(f))^2.
\]
Let $\{F_n\}_{n\geq 1}$ be a sequence of double Wigner integrals.
Then, as $n\to\infty$, we have 
\begin{enumerate}
\item[(i)] $F_n\overset{\rm law}{\to} A+I_2^S(f)$
\end{enumerate}
if and only if all the following are satisfied:
\begin{enumerate}
\item[(ii-a)] $\widehat{\kappa}_2(F_n)\to \widehat{\kappa}_2(A+I_2^S(f))=\mu_0^2+\|f\|^2_{L^2(\R_+^2)}$;
\item[(ii-b)] $\sum_{r=3}^{{\rm deg}Q} \frac{Q^{(r)}(0)}{r!}\,\widehat{\kappa}_r(F_n)
\to \sum_{r=3}^{{\rm deg}Q} \frac{Q^{(r)}(0)}{r!}\,\widehat{\kappa}_r(I_2^S(f))$;
\item[(ii-c)] $\widehat{\kappa}_r(F_n)\to \widehat{\kappa}_r(I_2^W(f))$ for  $a(f)$ consecutive values of $r$, with $r\geq 2(1+{\bf 1}_{\{\mu_0\neq 0\}})$.
\end{enumerate}
\end{theorem}

To conclude, we describe some easy consequences of Theorem \ref{moments-intdouble2}:
  \begin{enumerate}
  \item Consider first the situation where $\|f\|_{L^2(\R_+^2)}=0$ (that is, ${\rm rank}(f)=a(f)=0$) and $\mu_0\neq 0$. In this case, $Q(x)=x^4$ and condition $(ii-c)$ is
immaterial. Therefore, $F_n\overset{\rm law}{\to}
\mathcal{S}(0,\mu_0^2)$ if and only if $\widehat{\kappa}_2(F_n)\to \mu_0^2$ (condition $(ii-a)$) and $\widehat{\kappa}_4(F_n)\to 0$ (conditions $(ii-b)$). As such, one recovers a result by Kemp, Nourdin, Peccati and Speicher (in the case of double integrals), see \cite{knps}.
\item
Consider now the situation where, in Theorem \ref{moments-intdouble}, one has $r={\rm rank}(f)<\infty$, $a(f)=1$, $\lambda_1(f)=1$ and $\mu_0=0$. This corresponds to the case where the limit $I_2^S(f)$
is a centered free Poisson random variable with rate $r$. In this case, $Q(x)=x^4-2x^3+x^2$.
Therefore,  $F_n\overset{\rm law}{\to}
I_2^S(f)$ if and only if $\widehat{\kappa}_2(F_n)\to r$ (conditions $(ii-a)$ and $(ii-c)$) and
$\widehat{\kappa}_4(F_n)-2\,\widehat{\kappa}_3(F_n)\to -r$ (condition $(ii-b)$). As such, one recovers a theorem of Nourdin and Peccati in the special case of double integrals, see \cite{NP-poisson}.
\item
Finally, consider the situation $\mu_0=0$ and
$f=e_1-e_2$, where $e_1,e_2\in L^2(\R_+)$ are orthogonal and have norm 1
(thus, $r={\rm rank}(f)=a(f)=2$). This corresponds
to the case
where the limit is $S_1^2-S_2^2$,
with $S_1,S_2\sim\mathcal{S}(0,1)$ independent. That is, the limit is distributed according to the tetilla law.
In this case, $Q(x)=x^6-2x^4+x^2$.
Therefore,  $F_n\overset{\rm law}{\to} S_1^2-S_2^2$
if and only if $\kappa_2(F_n)\to 2$ (conditions $(ii-a)$ and $(ii-c)$), $\kappa_4(F_n)\to 2$ (condition $(ii-c)$) and $\kappa_6(F_n)-2\kappa_4(F_n)\to -2$ (condition $(ii-b)$). As such, one recovers a theorem of Deya and Nourdin in the special case of double integrals, see \cite{tetilla}.
\end{enumerate}

\bigskip

\noindent
{\bf Acknowledgments}. We are grateful to Giovanni Peccati for bringing to our attention reference \cite{russe} and for pointing out an error in a previous version (which  is  unfortunately the published version \cite{NP12},  see also \cite{NP12errata}) of this work.

\end{document}